%% file: is-5.tex
\documentclass[graybox]{svmult}
\usepackage{amsmath,amssymb,mathrsfs,url,amsfonts} 

\usepackage{mathptmx}       
\usepackage{helvet}         
\usepackage{courier}        
\usepackage{type1cm}        
\usepackage{makeidx}         
\usepackage{graphicx}        
\usepackage{multicol}        
\usepackage[bottom]{footmisc}

\makeindex             

\newcommand\DN{\newcommand}

\input{is-def.tex}

\titlerunning{Infinite-dimensional Stochastic Differential Equations with Symmetry} 

\title*{ 
Infinite-dimensional Stochastic Differential Equations with Symmetry. }
\author{Hirofumi Osada }
\institute{Department mathematics, Kyushu University \at Nishi-ku, Fukuoka 819-0395, Japan, \email{osada@math.kyushu-u.ac.jp}
}

\begin{document} \maketitle

\abstract{We review recent progress in the study of infinite-dimensional stochastic differential equations with symmetry. This paper contains examples arising from random matrix theory. \\
AMC2010: 60H110, 60J60, 60K35, 60B20, 15B52 }

\section{Introduction} \label{s:1}
We consider $ (\R ^d )^{\N }$-valued infinite-dimensional stochastic differential equations (ISDEs) of $ \mathbf{X}=(X^i)_{i\in\N }$ of the form

\begin{align}& \label{:1a}
dX_t^i = 
\sigma (X_t^i,\mathbf{X}_t^{i \diamondsuit })dB_t^i + 
b (X_t^i,\mathbf{X}_t^{i \diamondsuit })dt 
\quad (i\in\mathbb{N})
.\end{align}
Here $ \mathbf{B}=(B^i)_{i\in\N }$ is $ (\Rd )^{\N }$-valued standard Brownian motion. 
For  $ \mathbf{X}=(X^i)_{i\in\N }$ we set 
$ \mathbf{X}^{i \diamondsuit } = (X^j)_{j\in \N \backslash \{ i \} }$. 
Coefficients $ \sigma (x,\mathbf{y}) $ and $ b (x,\mathbf{y})$ 
are defined on a subset $ \SSw $ 
of $ \Rd \times (\Rd )^{\N }$ independent of $ i \in \N $. 
By definition $ \sigma (x,\mathbf{y}) $ is $ \R ^{d^2}$-valued, 
and $ b (x,\mathbf{y})$ is $ \Rd $-valued. We assume that 
$ \sigma (x,\mathbf{y}) $ and $ b (x,\mathbf{y})$ 
are symmetric in $ \mathbf{y}=(y_i)_{i\in\N }$ for each $ x \in \Rd $. 
Therefore, the set of \eqref{:1a} is referred to as 
\lq\lq ISDEs with symmetry''. 
In the present article, we review recent results in this regard. 
Using a Dirichlet form technique and an analysis on tail $ \sigma $-fields of configuration spaces, 
we prove the existence and pathwise uniqueness of strong solutions of the ISDEs of \eqref{:1a}. 
We emphasize that the coefficients are defined only on a {\em thin} subset in $ (\Rd )^{\N }$ and the state space of the solution $ \mathbf{X}$ is in this subset. Solving the ISDEs of \eqref{:1a} includes identifying such a subset.

%
%
%
%

Let $ \mathsf{S} = \{ \mathsf{s}=\sum_i  \delta_{s_i};\ 
\mathsf{s}(K) < \infty \text{ for any compact }K \subset \Rd 
\} $ be the configuration space over $ \Rd $. 
 $ \mathsf{S}$ is a Polish space equipped with the vague topology. 
With the symmetry of $ \sigma (x,\mathbf{y}) $ and $ b (x,\mathbf{y})$ in $ \mathbf{y}$, we regard $ \sigma $ and $ b $ as functions on $ \Rd \times \mathsf{S}$. 
%
%
%
We denote these by the same symbol such that 
$  \sigma (x,\mathbf{y}) = \sigma (x,\mathsf{y}) $ and 
$  b (x,\mathbf{y}) = b (x,\mathsf{y}) $, where 
$ \mathsf{y}= \sum_i \delta_{y_i}$ for $ \mathbf{y}=(y_i)$. 
%
Then we rewrite the ISDEs of \eqref{:1a} for $ \mathbf{X}=(X^i)_{i\in\N }$ as: 
\begin{align}\label{:isde}& 
\tcb{dX_t^i = \sigma (X_t^i, \mathsf{X}_t^{i\diamondsuit}) dB_t^i + 
b (X_t^i, \mathsf{X}_t^{i\diamondsuit}) dt \quad (i \in \N )}
.\end{align}
Here $ \mathsf{X}^{i\diamondsuit} = 
\sum_{j\not= i}^{\infty} \delta _{X^j} $, which 
is the $ \mathsf{S}$-valued process 
$\{ \sum_{j\not= i}^{\infty} \delta _{X_t^j} \}$. 
 $ \mathbf{X}$ is called the {\em labeled} dynamics, and the associated {\em unlabeled} dynamics $ \mathsf{X}$ is given by 
$ \mathsf{X} = \sum_{i\in\N } \delta _{X^i}$.

If $ \sigma $ is a unit matrix and $ b $ is 
given by a pair interaction 
$ \Psi (x-y)$, \eqref{:1a} becomes 
\begin{align}& \label{:Lang}
dX^i_t = dB^i_t + \frac{\beta}{2} \sum ^{\infty}_{j=1,j\ne i} 
\nabla \Psi (X ^i_t-X^j_t) dt 
. \end{align}
Here $ \beta > 0 $ is inverse temperature. 
Having $ \Psi $ of Ruelle class and $ \Psi \in C_0^3 (\Rd )$, Lang \cite{lang.1,lang.2} then solved \eqref{:Lang}, Fritz \cite{Fr} constructed non-equilibrium solutions for $ d \le 4 $, and Tanemura \cite{tane.2} provided solutions for hard core Brownian balls. 
The stochastic dynamics $ \mathbf{X} $ given by the solution of \eqref{:Lang} are called the interacting Brownian motions (IBM).

%
%
%
These solutions are strong solutions in the sense that $ \mathbf{X}$ are functionals of the given Brownian motions $\mathbf{B}$ and initial starting points $ \mathbf{s}$. The method used in these studies are based on the classic It$ \hat{\mathrm{o}}$ scheme. Hence, if $ \Psi $ is of long range such as a polynomial decay, then it is difficult to apply this scheme. Tsai \cite{tsai.14} solved \eqref{:Lang} for the Dyson model. He used very cleverly a specific monotonicity of the logarithmic potential and its one-dimensional structure. As for the weak solution, we present a robust method based on the Dirichlet form technique from \cite{o.isde}. 
We present a general theory to give $ \mu $-pathwise unique strong solutions applicable to the logarithmic interaction from \cite{o-t.tail}.

Thus, our demonstration is divided into two steps. In the first step, we obtain weak solutions of ISDEs \eqref{:1a}. That is, we construct solutions $ (\mathbf{X},\mathbf{B})$ satisfying \eqref{:1a} (see \sref{s:2}--\sref{s:4}). 
%
In the second step, we prove the existence of strong solutions and the $ \mu $-pathwise uniqueness. For this, we perform a fine analysis of the tail $ \sigma $-field of $ \mathsf{S}$ (see \sref{s:5} and \sref{s:6}). In \sref{s:7}, we give ISDEs arising from random matrix theory. 
In \sref{s:8}, we present the algebraic construction of the dynamics, and the coincidence of the algebraic dynamics with solutions of ISDEs. 
%
%
%

\section{Unlabeled dynamics: quasi-Gibbs property}\label{s:2}
We next construct a natural $ \mu $-reversible unlabeled diffusion, where $ \mu $ is a point process. The key point is the quasi-Gibbs property of $ \mu $, which we proceed to describe. 
%
 
Let $ \Sr = \{ x \in \Rd ; |x| < r  \} $. Let 
$ \map{\pi_r, \pi_r^c}{\mathsf{S}}{\mathsf{S}}$ be projections such that 
$\pi_r (\mathsf{s})=\mathsf{s}(\cdot \cap S_r)$, 
$ \pi_r^c(\mathsf{s})=\mathsf{s}(\cdot \cap S_r^c)$. 
For a point process $ \mu $, we set 
\begin{align*}&\tcb{
\mu _{r,\mathsf{t} }^m (\cdot )= 
\mu ( \pi _r (\mathsf{s})\in \cdot  |\mathsf{s}(S_r)= m ,\pi_r^c(\mathsf{s})= \pi_r^c(\mathsf{t})) 
}\end{align*}
Let 
$ \map{\Phi }{\Rd }{\mathbb{R}\cup\{ \infty \} }$ and 
$ \map{\Psi }{(\Rd )^2}{\mathbb{R}\cup\{ \infty \} }$ be potentials. We set 
$$ \mathcal{H}_r = 
\sum_{{\small s_i \in S_r }} \Phi (s_i) + 
 \sum_{{\small s_i,s_j \in S_r,  i<j}} \Psi (s_i, s_j)
 .$$
A point process $ \mu $ is called a canonical Gibbs measure if 
$ \mu $ satisfies Dobrushin-Lanford-Ruelle (DLR) equation, that is, 
for $ \mu $-a.s.\! $ \mathsf{t}=\sum_j \delta _{t_j} $ 
\begin{align}& \label{:DLR}
\mu _{r,\mathsf{t} }^m = 
{c^m_{r,\mathsf{t} }}  e^{-\mathcal{H}_r - 
\sum_{x_i\in S_{r}, t _j \in S_r^c} \Psi (x_i,t _j)} 
d\Lambda ^m_r 
.\end{align}
Here $ \Lambda ^m_r = \Lambda ( \cdot | \mathsf{s} (S_r ) = m )$ and 
$ \Lambda _r $ is the Poisson PP with intensity $ 1_{S_r}dx $.

Point processes appear in random matrix theory in the form sine, Airy, Bessel, and Ginibre point processes having logarithmic potentials 
$$ \Psi (x,y)= - \beta \log |x-y| .$$
However, the DLR equation \eqref{:DLR} does not make sense for a logarithmic potential. 
Hence we introduce the notion of quasi-Gibbs measures: 
\begin{definition}\label{dfn:1}
$ \mu $ is $ (\Phi ,\Psi )$-quasi-Gibbs measure if $\exists \ c^m_{r,\mathsf{t} }$ such that 
\begin{align*} & \tcb{ 
{ c^m_{r,\mathsf{t} }} ^{-1}  e^{-\mathcal{H}_r } d\Lambda ^m_r \le 
\mu_{r,\mathsf{t} }^m \le {c^m_{r,\mathsf{t} }} e^{-\mathcal{H}_r } d\Lambda ^m_r }
\end{align*}
\end{definition}
By definition a canonical Gibbs measure is a quasi-Gibbs measure. 
We refer to \cite{o.rm,o.rm2} for a sufficient condition for quasi-Gibbs property. 
We assume: \\
\As{A1}  $ \mu $ is a quasi-Gibbs measure with upper semi-continuous 
$ (\Phi ,\Psi )$. Furthermore, 
$  a (x,\mathsf{s}) = \sigma ^t (x,\mathsf{s}) \sigma (x,\mathsf{s}) $ 
 is bounded and uniformly elliptic. 
\\
\As{A2} There exists a $ 1 < p \le \infty $ such that 
the $ k $-point correlation function $ \rho ^k $ of $ \mu $ 
is in $ L_{\mathrm{loc}}^p ((\Rd )^m)$ for each $ k \in \N $.

For a given point process $ \mu $ we introduce 
 a Dirichlet form such that 
\begin{align}& \label{:2b}
\mathcal{E}^{\mu } (f,g)= \int_{\mathsf{S}} 
\mathbb{D}[f,g]d\mu  
,\quad 
\mathbb{D}[f,g] = 
\frac{1}{2} \sum_{i} 
a (s_i,\mathsf{s}^{i\diamondsuit })\PD{\check{f}}{s_i}\cdot \PD{\check{g}}{s_i}
.\end{align}
Here we set $ \mathsf{s}^{i\diamondsuit }=\sum_{j\not=i} \delta_{s_j}$ 
for $ \mathsf{s}=\sum_i\delta_{s_i}$, 
and 
$ f (\mathsf{s}) = \check{f} (s_1,s_2,\ldots )$, where 
$ \check{f}$ is symmetric in $ (s_1,s_2,\ldots )$. 
Note that $ \mathbb{D}[f,g] $ is a function of $ \mathsf{s}$ by construction. 
\begin{theorem}[\cite{o.dfa,o.rm,o-t.tail}]	\label{l:diri} 
Let $ \dcirc $ be the set of local, smooth functions on $ \mathsf{S}$. 
Set $ \dcirc ^{\mu} = \{ f \in \dcirc \cap L^2(\mu)\, ;\, \mathcal{E}^{\mu } (f,f) < \infty   \} $. 
\\
\thetag{i}  Assume \As{A1}. Then 
$ ( \mathcal{E}^{\mu} ,\dcirc ^{\mu} )$ is closable on $ L^2(\mu)$. 
\\
\thetag{ii} Assume \As{A1} and \As{A2}. Then there exists 
a diffusion $ \mathsf{X}_t=\sum_i \delta_{X_t^i}$ associated with the closure 
$ ( \mathcal{E}^{\mu} ,\mathcal{D}^{\mu} ) $ of 
$  ( \mathcal{E}^{\mu} ,\dcirc ^{\mu} )$ on $ L^2(\mu)$. 
\end{theorem}
The local boundedness of the correlation functions is used for the quasi-regularity of the Dirichlet form. Once quasi-regularity is established, the existence of $ \mu $-reversible diffusion is immediate from the general theory \cite{m-r,FOT.2}. 

Unlabeled dynamics are also obtained in \cite{akr,y.96} with a different frame work. It is now proved these are the same dynamics as in \cite{o-t.sm,o-t.core}. We remark that ergodicity of unlabeled dynamics with grand canonical Gibbs measures with small enough activity constant 
is obtained in \cite{akr}.

\section{Labeled dynamics: A scheme of Dirichlet spaces}
\label{s:3}
We next lift the unlabeled dynamics $ \mathsf{X}$ in \tref{l:diri} to a labeled dynamics $ \mathbf{X}=(X^i)_{i\in\N}$ solving \eqref{:1a}. 
For this we present a natural scheme of Dirichlet spaces describing the labeled dynamics $ \mathbf{X}$. 
We assume a pair of mild assumptions: 

\noindent 
\As{A3} $ \{ X^i \} $ do not collide with each other 
  \tcb{(non-collision)}
\\
\tcb{\As{A4}} each tagged particle $ X^i $ never explode 
 \tcb{(non-explosion)} 

Let 
$ \mathsf{S}_{s,i}=\{ \mathsf{s} \in \mathsf{S}\, ;\, 
\mathsf{s}(\{ x \} )=0  \text{ for all } x \in S ,\ 
\, \mathsf{s}( S ) = \infty \}  $. Then \As{A3} is equivalent to 
$ \mathrm{Cap}^{\mu } (\mathsf{S}_{s,i}^c) = 0 $. 
\As{A4} follows from $ \rho ^1(x) = O (e^{|x|^{\alpha }})$, 
$ \alpha < 2$. 

We call $ \ulab  $ the unlabeling map if $ \ulab ((s_i)) = \sum_i \delta _{s_i} $. We call $ \lab $ a label if $ \lab $ 
is defined for $ \mu $-a.s.\, 
$ \mathsf{s}$, and $ \ulab \circ \lab (\mathsf{s}) = \mathsf{s}$. 
For a unlabeled dynamics satisfying \As{A3} and \As{A4}, 
the particles can keep the initial label $ \lab (\mathsf{s})$. 
Thus 
we can construct a map $ \lpath $ to $ C([0,\infty);\RdN )$ 
such that 
$ \{ \ulab (\lpath (\mathsf{X}_t))\}_{t\in[0,\infty )} = \mathsf{X}$. 
Hence we obtain: 
\begin{theorem}[{\cite{o.tp}}] \label{l:31}
Assume \As{A1}--\As{A4}. 
Then there exists a labeled dynamics 
$ \mathbf{X}  = (X ^i)_{i\in\N }$ such that 
$ \mathsf{X}  = \sum_{i\in\mathbb{N}} \delta_{X ^i}$ and that 
$ \mathbf{X}_0 = \lab (\mathsf{X}_0)$. 
\end{theorem}

Remark that $ \RdN $ has no good measures. Then 
no Dirichlet forms on $ \RdN $ associated with the labeled dynamics $ \mathbf{X}$. 
We hence introduce the scheme of spaces $ (\Rd )^m \times \mathsf{S}$ 
with Campbell measures $ \mu ^{[\M ]}$ such that 
$d\mu ^{[\M ]} = \rho ^{\M }(\mathbf{x}_{\M }) 
\mu _{\mathbf{x}_{\M }} (d\mathsf{s})d\mathbf{x}_{\M } $, 
where 
$  \rho ^{\M }$ is a $\M $-point correlation function of $ \mu $ and 
$ \mu _{\mathbf{x}_m}$ is the reduced Palm measure conditioned at $ \mathbf{x}_{\M }$. 
For $ a = \sigma^t\sigma $, let $ \mathbb{D}^{[\M ]}$ be the square field on $\Rdm  \times \mathsf{S}$ defined similarly as $ \mathbb{D}$ on $ \mathsf{S}$ given by \eqref{:2b} in \sref{s:2}.  Let 
\begin{align*}&
\mathcal{E}^{[\M ]} (f,g)= 
\int_{\Rdm  \times \mathsf{S}} \mathbb{D}^{[\M ]} [f,g] d\mu ^{[\M ]} 
.\end{align*}
Let $ \dcirc ^{[1],\mu }= \{ f \in C_0^{\infty}(\Rd )\ot \dcirc ; 
\mathcal{E}^{[\M ]} (f,f) < \infty ,\, f \in L^2 (\mu ^{[\M ]})
\} $. 
\begin{theorem}[{\cite{o.tp}}] \label{l:32}
Assume \As{A1} and \As{A2}. 
Then $ (\mathcal{E}^{[\M ]} ,\dcirc ^{[1],\mu })$ is 
closable on $ L^2 (\mu ^{[\M ]} )$, and  
its closures is quasi-regular.
Hence the associated diffusion 
$ (\mathbf{X} ^{[\M ]}, \mathsf{X} ^{[\M ]*})$ exists. Here we write 
$ (\mathbf{X} ^{[\M ]}, \mathsf{X} ^{[\M ]*}) = 
(X ^{[\M ],1},\ldots,X ^{[\M ],\M }, 
\sum_{i= \M +1}^{\infty} \delta_{X ^{[\M ],i}})$. 
\end{theorem}

Let $ (\mathcal{E}^{\mu}, \mathcal{D}^{\mu } , L^2(\mu) )$ be 
the original Dirichlet form. Let 
$ \mathsf{X}= \sum_{i=1}^{\infty} \delta_{X^i}$ be the associated unlabeled diffusion. 
We fix a label $ \lab $. 
Let $ \mathbf{X}=(X^i)_{i\in \N}$ be the labeled dynamics given by $ \lab $. 
We set $ (\mathbf{X}^m,\mathsf{X}^{m*}) = (X^{1 },\ldots, X^{\M }, \sum_{i =\M +1}^{\infty} \delta_{X^{i}}) $. 
%
%
\begin{theorem}[\cite{o.tp}]\label{l:coup2}
Assume \As{A1}--\As{A4}. 
Assume $ (\mathbf{X}^{[\M ]}, \mathsf{X}^{[\M ]*}) $ and 
$ (\mathbf{X}^{\M }, \mathsf{X}^{\M *}) $ start at the same initial point. Then 
$ (\mathbf{X}^{[\M ]}, \mathsf{X}^{[\M ]*}) = 
(\mathbf{X}^{\M }, \mathsf{X}^{\M *}) $ 
in distribution for each $  \M \in \N $. 
\end{theorem}
Instead of the huge space $ \RdN $, 
we use a scheme of countably infinite {\em good} infinite-dimensional spaces 
$\{  (\Rd )^{m}\times \mathsf{S}\}_{m\in  \{ 0 \} \cup \N }$. 
Using the diffusion $ \mathsf{X}$ on the original unlabeled space $ \mathsf{S}$, we construct a scheme of the {\em coupled} diffusions $ (\mathbf{X}^{\M }, \mathsf{X}^{\M *})$ on $ \Rdm \times \mathsf{S}$ associated with the scheme of Dirichlet spaces 
$ ( \mathcal{E}^{[\M ]}, L^2(\mu ^{[\M ]})) \text{ on } (\Rd )^{\M }\times \mathsf{S} $. 
This construction is key for the ISDE-representation below.

\section{ISDE-representation: Logarithmic derivative}\label{s:4}
\begin{definition}[\cite{o.isde}]\label{dfn:41}
Let $ \nabla _x  $ be the nabla on $ \mathbb{R}^d$. 
 $ \mathsf{d}^{\mu } \in 
L_{\mathrm{loc}}^1( \R ^d \times \mathsf{S} , \mu ^{[1]} )^d $ is 
called the {\em logarithmic derivative} of $ \mu $ if, 
for all 
$ f \in C_0^{\infty}(\R ^d)\ot L^{\infty}(\mathsf{S}) $,  
\begin{align}& \label{:LD}
\int _{\R ^d \times \mathsf{S} }\nabla _x f d\mu ^{[1]} = 
 - \int _{\R ^d \times \mathsf{S} } f  \mathsf{d}^{\mu }  d\mu ^{[1]} 
.\end{align}
\end{definition}
%

Let 
$  a (x,\mathsf{s}) = \sigma ^t (x,\mathsf{s}) \sigma (x,\mathsf{s}) $ 
as before. We set $ \nabla_x a $ such that 
$ \nabla_x a = \sum_{j=1}^d \PD{}{x_j}a_{ij}$. 
We introduce a \lq\lq geometric'' differential equation on 
$ \mathsf{d}^{\mu } (x,\mathsf{s})
 =: \nabla_x \log \mu ^{[1]} (x,\mathsf{s})$:
\begin{align}&\label{:GDE}
\nabla_x a(x,\mathsf{s}) + 
a(x,\mathsf{s}) \nabla_x \log \mu ^{[1]} (x,\mathsf{s}) 
 = 
2 b (x, \mathsf{s}) 
.\end{align}
\As{A5} $ \mu $ has a logarithmic derivative $ \mathsf{d}^{\mu }$. \\
\As{A6} The logarithmic derivative $ \mathsf{d}^{\mu }$ satisfies 
\eqref{:GDE}. 
\begin{theorem}[\cite{o.isde}]	\label{l:isde} 
Assume \As{A1}--\As{A6}. 
Then there exists an $ \mathsf{S}_0 \subset \mathsf{S}$ 
such that $ \mu (\mathsf{S}_0) = 1$ 
and that, for each  $ \mathbf{s}\in \ulab ^{-1}(\mathsf{S}_0) $,   ISDE \eqref{:1a}  has 
a solution $ (\mathbf{X},\mathbf{B})$ satisfying $ \mathbf{X}_0 = \mathbf{s} $ and 
$ \mathbf{X}_t \in \ulab ^{-1}(\mathsf{S}_0) $  for all $ t $. 
\end{theorem}

From the coupling in \tref{l:coup2} and Fukushima decomposition (It$ \hat{\mathrm{o}}$ formula), we prove that 
$ \mathbf{X}=(X ^i)_{i\in\N }$ satisfies the ISDEs of \eqref{:1a}. 
We use the $ \M $-labeled process 
$ (\mathbf{X}^{\M }, \mathsf{X}^{\M *}) $, 
to apply  It$ \hat{\mathrm{o}}$ formula to coordinate functions 
$ x_1,\ldots,x_{\M }$.

\section{Strong solutions of ISDEs and pathwise uniqueness}\label{s:5}
We lift the weak solutions $ (\mathbf{X},\mathbf{B}) $ to pathwise unique strong solutions. 
For this purpose, we introduce a scheme consisting of 
{\em an infinite system of finite-dimensional SDEs with consistency} (IFC),  and 
perform an analysis of the tail $ \sigma $-field of the path space 
$  \WTN = C([0,T]; \SN )$. 
The key idea is the following interpretations: 

\noindent \quad \quad \bub 
\tcb{a single ISDE} $ \Longleftrightarrow $ 
a scheme of \tcb{IFC}. 

\noindent \quad \quad \bub 
the tail $ \sigma $-field of $ \WTN $ $ \Longleftrightarrow $ 
the boundary condition of the ISDEs. 
\\
The method is robust and may be applied to many other models. 
We consider {\em non-Markovian} ISDEs because the argument is general.  
Let $ (\WTsol ,\Sz , \{ \sigma ^i \}, \{ b^i \} )$ be 
\begin{align*}
&  \text{$ \WTsol $F a Borel subset of $ \WTN $}
&& \text{\tcr{(space of solutions of the ISDEs)}}
,\\
&  \text{$ \Sz  $: a Borel subset of $ \SN $}
&& \text{\tcr{(initial starting points of the ISDEs)}}
,\\
&  \map{\sigma ^i , b^i }{\WTsol  }{\WTN }
&& \text{\tcr{(coefficients of the ISDEs)}}
.\end{align*}
We consider the ISDEs on $ \SN  $ of the form: 
$ \mathbf{X} = \{ (X_t^i )_{i\in\mathbb{N}} \}_{t \in [0,T]} \in \WTN $ 
\begin{align}\label{:5a} & 
dX_t^i = \sigma ^{i}(\mathbf{X})_t dB_t^i + 
b^{i}(\mathbf{X})_t dt \quad (i\in\mathbb{N})
,\\ \notag  
&
\mathbf{X}_0=\mathbf{s} =(s_i)_{i\in\mathbb{N}} \in  \Sz  
,\quad \mathbf{X} \in \WTsol 
.\end{align}
Note that \eqref{:1a} is a special case of \eqref{:5a}. We assume: 

\smallskip
\noindent 
\As{P1} 
The ISDEs \eqref{:5a} has a weak solution $ (\mathbf{X},\mathbf{B}) $. 
\quad (not a strong solution!) 

\smallskip

From a weak solution $ (\mathbf{X},\mathbf{B})$, we define a new SDE of 
$ \mathbf{Y}^{m} = (Y_{t}^{m,i})_{i=1}^m $ 
such that 
\begin{align}\label{:5b}&
dY_t^{m,i} = \sigma ^{i} (\mathbf{Y}^{m} , \mathbf{X}^{m*})_t dB_t^i + 
b^{i}(\mathbf{Y}^{m} , \mathbf{X}^{m*})_t dt 
\quad  ( i=1,\ldots,m) 
,\\ \notag &
\mathbf{Y}_0^{m} = (s_1,\ldots,s_m) 
,\quad (\mathbf{Y}^{m} , \mathbf{X}^{m*}) \in \WTsol 
\end{align}
for each $ \mathbf{X} \in \WTsols $, 
$ \mathbf{s}=(s_i)_{i=1}^{\infty} \in \Sz  $, and $ m \in \mathbb{N}$. 
Here $ \mathbf{X}^{m*}:=(X^{n})_{n>m}$ 
is interpreted as a part of the coefficients of the SDE \eqref{:5b} and 
$ \WTsols  = \{ \mathbf{X}\in \WTsol ; 
\mathbf{X}_{0} = \mathbf{s} \} $. 
Indeed, we regard \eqref{:5b} as {\em finite-dimensional}
 SDEs of $ \mathbf{Y}^{m}$. \eqref{:5b} become automatically time-inhomogeneous SDEs. We have therefore obtained a scheme of finite-dimensional SDEs of $ \{ \mathbf{Y }^m \}_{m\in\N } $. 
We assume: 

\smallskip 

\noindent 
\As{P2} The SDE \eqref{:5b} has a unique, strong solution 
for each $ \mathbf{s} \in \Sz  $, 
$ \mathbf{X} \in \WTsols $, and $ m \in \mathbb{N} $. 

\smallskip 

Let  $ \Pbars $ be the distribution of solution 
$ (\mathbf{X},\mathbf{B})$ on $ \WW $. 
Let 
\begin{align*}
&\text{
$ \PbarsB  = \Pbars (\mathbf{X} \in \cdot | \mathbf{B}) $ 
 and 
 $\PBr = \Pbars (\mathbf{B} \in \cdot ) $, }
\\&
\text{$ \TTT  = \cap_{m=1}^{\infty} \sigma [\mathbf{X}^{m*}]  $}
,\\&
\text
{$ \Tone [\PsB ] = \{ \mathbf{A}\in \TTT \, ;\, \PsB (\mathbf{A}) = 1 
\} $. }
\end{align*}
\As{P3} $ \TTT  $ is $ \PsB $-trivial for each $ \mathbf{s} \in \Sz  $ 
and 
$ \PBr $-a.s.\! $ \mathbf{B}$. 
\begin{theorem}[\cite{o-t.tail}] \label{l:q1} 
Assume \As{P1}--\As{P3}. 
Then \\
\thetag{i} 
$ \mathbf{X}$ is a strong solution of the ISDEs \eqref{:5a} 
for each $ \mathbf{s}\in \mathbf{S}_{0}$. 
\\
\thetag{ii} 
Let $ \mathbf{X}_{\mathbf{s}}$ and $ \mathbf{X}_{\mathbf{s}}' $ 
be strong solutions of the ISDEs \eqref{:5a} starting at 
$ \mathbf{s}\in\mathbf{S}_{0}$ 
 defined on the same space of Brownian motions $ \mathbf{B}$. 
Then $ \mathbf{X}_{\mathbf{s}} =  \mathbf{X}_{\mathbf{s}}'  $ for 
$ \PBr $-a.s.\! $ \mathbf{B}$ 
 if and only if  $\Tone [\PsB ] = \Tone [\PsB ']$ for 
$ \PBr $-a.s.\! $ \mathbf{B}$. 
\end{theorem}

\noindent {\sf 
Idea of the proof \thetag{i}: }  
Let 
$ (\mathbf{X},\mathbf{B})$ be a weak solution of ISDE given by \As{P1}, 
and fix it. 
Let 
$ \mathbf{Y}^{m}$ be the unique strong solution of \eqref{:5b} 
given by \As{P2}. 
By construction 
$ \mathbf{Y}^{m}$ is 
$ \sigma[\mathbf{B}]\bigvee \sigma[\mathbf{X}^{m*}]$-measurable. 
Because the solution \eqref{:5b} is unique, we see that 
$ \mathbf{Y}^{m} = \mathbf{X}^{m} :=(X^1,\ldots,X^m)$. 
Let $ \mathbf{Y}$ be the limit $ \mathbf{Y}=\limi{m} \mathbf{Y}^{m}$. 
Then $ \mathbf{Y}=\mathbf{X} $ and $ \mathbf{Y}$ is 
$ \sigma[\mathbf{B}]\bigvee \TTT  $-measurable. 
Because $ \TTT $ is $ \PsB $-trivial by \As{P3}, 
 $ \mathbf{Y}$ depends only on $ \mathbf{s}$ and $ \mathbf{B}$. 
This means $ \mathbf{Y}=\mathbf{X}$ is a strong solution. \qed

In \cite{o-t.tail} we intoduce a notion of IFC solution, with which 
we generalize \tref{l:q1}. 

\section{Tail triviality: Application to interacting Brownian motions.}
\label{s:6}
We return to the Markovian-type ISDEs of \eqref{:1a}. 
We assume \As{A1}--\As{A6}. 
We apply \tref{l:q1} to ISDEs of \eqref{:1a} 
by checking \As{P1}--\As{P3}. 
\As{P1} follows from \tref{l:isde}. 
Controlling the capacity of 
$ (\mathcal{E}^{\mu}, \mathcal{D}^{\mu } , L^2(\mu) )$ 
we obtain \As{P2}. 
Because $ \mathbf{X} \in \WTsol $ and 
$ \WTsol $ is a nice subset of $ \WT (\RdN ) $, 
we can assume \As{P2} for the solution of \eqref{:1a}. 
Dirichlet form theory proves that $ \mathbf{X}$ stays in $ \WTsol $. 
Indeed, such a condition is reduced to a calculation of capacity 
related to the unlabeled Dirichlet space \cite{FOT.2,m-r}. 
Roughly speaking,  \As{P2} is satisfied if 
$ \nabla _{x}^m \mathsf{d}^{\mu } 
\in \mathcal{D}_{\mathrm{loc}}^{\mu ^{[1]}} $ 
for a suitable $ m $, see  \cite[Sections 8,9]{o-t.tail}. 
\begin{theorem}[\cite{o-t.tail}]	\label{l:tail3} 
Assume \As{Q1}--\As{Q3} below. 
 Then \As{P3} holds. 

\noindent 
\As{Q1} $ \mu $ is tail trivial. That is, 
$ \mu (A) \in \{ 0,1 \} $ for all 
$ A \in \mathcal{T} (\mathsf{S} ):=\cap_{r=1}^{\infty} \sigma [\pi_r^c]$. 
\\
\As{Q2} 
$ P_{\mu}\circ \mathsf{X}_t^{-1} \prec \mu $ for all $ t$. 
\quad (absolute continuity condition). 
\\
\As{Q3} $ P_{\mu } (\cap_{r=1}^{\infty} \{ \mathsf{m}_r (\mathsf{X}) < \infty \} ) = 1 $, \\
where 
$ \mathsf{m}_r = \inf\{ m\in\mathbb{N}; X^i \in C([0,T]; S_r^c) \text{ for } m < \forall i \in \mathbb{N}\}$ for $ \mathsf{X} = \sum_{i\in\mathbb{N}} \delta_{X^i}
$. 
\end{theorem}

\begin{remark}\label{r:6} 
\thetag{i} Determinantal point processes satisfy \As{Q1} (see \cite{o-o.tt}). 
\\\thetag{ii} 
\As{Q2} is obvious because the unlabeled dynamics is $ \mu $-reversible. 
\\\thetag{iii} 
\As{Q3} is satisfied if the one-point correlation function $ \rho ^1$ satisfies$ \rho ^1(x) = O (e^{|x|^{\alpha }}) $ $ (|x|\to \infty)$ 
 for some $ \alpha < 2 $. 
\end{remark}

Let $ \mathbf{T} = \{ \mathbf{t}=(t_1,\ldots,t_m)\, ;\, t_i \in [0,T], m \in \mathbb{N}\} $ and 
$ \mathbf{X}_{\mathbf{t}}^{n*} = 
(\mathbf{X}_{t_1}^{n*},\ldots,\mathbf{X}_{t_m}^{n*})$. 
Let 
\begin{align}\notag &
\5 (\mathsf{S}) = \bigvee_{\mathbf{t} \in \mathbf{T}}
\bigcap_{r=1}^{\infty}\sigma [\pi _r^c (\mathsf{X}_{\mathbf{t}}) ]
,\quad 
\5 (\SN ) = \bigvee_{\mathbf{t} \in \mathbf{T}}
\bigcap_{n=1}^{\infty}\sigma [\mathbf{X}_{\mathbf{t}}^{n*}]
.\end{align}
Hence by definition, 
$ \5 (\mathsf{S}) $ is the cylindrical tail $ \sigma $-field of the unlabeled path space 
and 
$ \5 (\SN ) $ is the cylindrical tail $ \sigma $-field of the labeled path space $ \WTSN $. 
%
We  deduce the triviality of $ \TpathTSN $ from that of $ \mathsf{S}$. 
We do this step-by-step following the scheme: 
\begin{align*}&
\mathcal{T}(\mathsf{S}) 
\xrightarrow[\As{Q1}, \,  \As{Q2}]{(\textrm{Step I})}
&& \5 (\mathsf{S})
\xrightarrow[\As{Q3}]{(\textrm{Step II})}
&& \5 (\SN ) 
\xrightarrow[\thetag{IFC}]{(\textrm{Step III})}
&& \TpathTSN  
\\ \notag &  {\mu }
&& \mathsf{P}_{\mu }
&& \mathbf{P}_{\mu^{\lab }}= \int \Pbars (\mathbf{X} \in \cdot ) d\mu ^{\lab }
&& \PsB \text{ a.s.\!\! } (\mathbf{s},\mathbf{B})  
.\end{align*}
%
We denote by 
$ \mathcal{B}(\SS )^{\mu } $ the completion of the $ \sigma $-field $ \mathcal{B}(\SS ) $ 
with respect to $ \mu $. %

\begin{definition}\label{d:21Z}
For $ \mathbf{s}\in \SN $, we set 
$ \mathbf{X}_{\mathbf{s}} = \{\mathbf{X}_{\mathbf{s},t}\}_{ t\in [0,\infty)} $ 
such that $ \mathbf{X}_{\mathbf{s},0} = \mathbf{s}$. 
We set  
$$ \{\mathbf{X}_{\mathbf{s}}\} _{\mathbf{s} \in \lab (\mathsf{H}) } 
= \{\{\mathbf{X}_{\mathbf{s},t}\}_{ t\in [0,\infty)} 
\} _{\mathbf{s}  \in \lab (\mathsf{H}) }.$$
\thetag{i} We call $\3 $ a $ \mu $-solution of \eqref{:1a} if 
$ \mathsf{H} \in  \mathcal{B}(\SS )^{\mu } $ satisfying 
$ \mu (\mathsf{H}) = 1 $ and $ \ulab ^{-1}(\mathsf{H}) \subset \SSw $, and if 
$ ( \mathbf{X}_{\mathbf{s}},\mathbf{B}) $ is a solution of 
 ISDE \eqref{:1a} for each $ \mathbf{s} \in \lab (\mathsf{H}) $. 
\\\thetag{ii} 
We call 
$\3 $  
a $ \mu $-strong solution if it is a $ \mu $-solution such that 
$ ( \mathbf{X}_{\mathbf{s}},\mathbf{B}) $ is a strong solution 
 for each $ \mathbf{s} \in \lab (\mathsf{H}) $. 
\end{definition}

\begin{definition}\label{d:21B}
We say that the $ \mu $-strong uniqueness holds if the following holds. \\
\thetag{i} 
The $ \mu $-uniqueness in law holds. 
That is, 
$ \mathbf{X}_{\mathbf{s}} = \mathbf{X}_{\mathbf{s}}' $ in law 
for each $ \mathbf{s} \in \lab (\mathsf{H}\cap\mathsf{H}')$ 
for any pair of $ \mu $-solutions  $\3 $ and 
$ (\{\mathbf{X}_{\mathbf{s}}'\}_{\mathbf{s}\in\lab (\mathsf{H}')},
 \mathbf{B}')$ 
satisfying \As{Q2}. 
\\
\thetag{ii} 
A $ \mu $-solution $\3 $ 
satisfying \As{Q2} is a $ \mu$-strong solution 
$ (\{\mathbf{X}_{\mathbf{s}}\}_{\mathbf{s}\in\lab (\mathsf{H}')}, 
\mathbf{B})$
for some $ \mathsf{H}' \subset \mathsf{H}$. 
\\
\thetag{iii} 
The $ \mu $-pathwise uniqueness holds. 
That is, 
$ P^{\mathbf{B}} (\mathbf{X}_{\mathbf{s}} = \mathbf{X}_{\mathbf{s}}') = 1 $ for each $ \mathbf{s}\in  \lab (\mathsf{H}\cap\mathsf{H}' )$, where 
$ \{\mathbf{X}_{\mathbf{s}}\}_{\mathbf{s}\in\lab (\mathsf{H})}$ and 
$ \{\mathbf{X}_{\mathbf{s}}'\}_{\mathbf{s}\in\lab (\mathsf{H}')}$ 
 are any pair of $ \mu $-strong solutions 
defined for the same Brownian motion $ \mathbf{B}$ satisfying \As{Q2}. 
\end{definition}

\begin{theorem}[{\cite{o-t.tail}}]	\label{l:8}
Make the same assumptions as for \tref{l:isde}. 
Assume \As{P2}, \As{Q1}, and \As{Q3}. 
Then \eqref{:1a} has a $ \mu $-strong solution $ \mathbf{X}$ 
such that the associated unlabeled dynamics $ \mathbf{X}$ 
is $ \mu $-reversibile, and the $ \mu $-strong uniqueness holds. 
\end{theorem}

\section{Examples arising from random matrix theory. }\label{s:7}
The first three examples are particle systems in $ \mathbb{R}$ ($ [0,\infty)$ for Bessel), whereas the last example is in $ \mathbb{R}^2$. All examples have logarithmic interaction potential. 

\noindent 
{\bf Sine, Airy, and Bessel IBM \cite{o.isde,tsai.14,o-t.airy,o-h.bes}: } 
Let  $ d=1$. 
\begin{align}& \tag{{Dyson model, Sine}}
\tcb{dX_t^i = dB_t^i + \frac{\beta }{2}\lim_{r\to\infty }
\sum_{|X_t^i-X_t^j |<r ,\  \scriptstyle j\not= i} \frac{1}{X_t^i-X_t^j}dt 
}
,\\ \tag{{Airy}} &
\tcb{ dX^i_t=dB^i_t 
+ \frac{\beta}{2} \lim_{r\to\infty} \{ ( 
 \sum_{j \not= i, \ |X^j_t |<r}\frac{1}{X^i_t -X^j_t } ) 
- \frac{1}{\pi }\int_{-r}^0\frac{\sqrt{-x}}{-x}dx 
 \}dt }
, \\
\tag{Bessel} & 
 dX_t^i = dB_t^i + 
\frac{a}{2X_t^i}dt 
+ 
\frac{\beta }{2}\sum_{j\not= i}^{\infty} 
\frac{1}{X_t^i-X_t^j}dt 
 \quad \text{  $ a \ge 1$}
.\end{align}
The equilibrium states of these dynamics are sine, Airy, and Bessel point processes ($ \beta = 1,2,4$ (sine, Airy), $\beta = 2 $ (Bessel)). 
 These point processes correspond to bulk, soft edge, and hard edge scaling limits respectively. 
The relationships to inverse temperature are: 
 $ \beta =1 \Rightarrow $GOE, $ \beta =2 \Rightarrow $GUE, and $ \beta =4 \Rightarrow $GSE, respectively.  

\noindent 
 {\bf Ginibre IBM \cite{o.isde}: } Let $ d=2 $ and $ \beta = 2 $. 
We consider two ISDEs. 
\begin{align}
& \label{:Gin1}
dX_t^i = dB_t^i + 
 \lim_{r\to\infty }
\sum_{|X_t^i - X_t^j|<r ,\ \scriptstyle j\not= i} 
\frac{X_t^i-X_t^j}{|X_t^i-X_t^j|^2}dt 
,\\
& \label{:Gin2} 
dX_t^i = dB_t^i  - X_t^i + 
 \lim_{r\to\infty }
\sum_{|X_t^j|<r ,\ \scriptstyle j\not= i} 
\frac{X_t^i-X_t^j}{|X_t^i-X_t^j|^2}dt 
.\end{align}

\noindent 
The equilibrium state $ \mug $ is the Ginibre point process, which has various rigidities such as small variance \cite{shirai.06}, number rigidity \cite{GP}, and dichotomy in its reduced Palm measures \cite{o-s.abs}. 
%
\noindent 
The drift coefficients are equal on and tangential to the support of $ \mug $, yielding the coincidence of the solutions of \eqref{:Gin1} and \eqref{:Gin2}. 
This dynamical rigidity reflects rigidity of $ \mug $. 

Another dynamical rigidity is the sub-diffusivity \cite{o.sub}: 
$ \limz{\epsilon} \epsilon X_{t/\epsilon ^2}^i =  0$.
All translation invariant IBMs in $ \Rd $ ($ d \ge 2 $) with Ruelle-class potentials with hard core are {\em diffusive} \cite{o.p}. 
Therefore the sub-diffusivity is in contrast with Ruelle-class potentials, 
and indicates a dynamical rigidity as a slow down in tagged particles.

\section{Algebraic construction and finite particle approxiations} \label{s:8}
An algebraic construction is known for 
stochastic dynamics related to point processes 
appearing in random matrix theory in $ \R $ 
with $ \beta = 2$, 
which is given by space-time correlation functions {\it e.g.} \cite{j.02,KT07b,KT11,KT11-b}. For example, 
%
as for the Airy$_2 $ point process, 
the multi-time, moment generating function is 
\begin{align}&\notag 
{\bf E} \big[\exp \big\{ \sum_{m=1}^{M} 
 \langle f_{t_m},  \mathsf{X}_{t_m}\rangle  \big\} \big]=
 \mathop{\rm Det}_
{\substack
{(s,t)\in \mathbf{t}^2, \, (x,y)\in \R^2}
}
\Big[\delta_{s t} \delta(x-y)
+ {\bf K}_{\mathrm{Airy}}(s, x; t, y) \chi_{t}(y) \Big]
.\end{align}
Here 
$ \mathbf{t}=\{ t_1,\ldots,t_M \} $, $ \chi_t = e^{f_t}-1$, and 
$ {\bf K}_{\Ai}$ is the extended Airy kernel 
\begin{align*}
&{\bf K}_{\Ai}(s,x;t,y) 
=\begin{cases}
\displaystyle{
\int_{0}^{\infty} du \, e^{-u(t-s)/2} \Ai(u+x) \Ai(u+y),
}
&t \geq s
\cr 
\displaystyle{
- \int_{-\infty}^{0} d u \, e^{-u(t-s)/2} \Ai(u+x) \Ai(u+y),
} 
&t < s.
\end{cases}
\end{align*}
\begin{theorem}[\cite{o-t.sm,o-t.core}]	\label{l:alg} 
The algebraic construction and the ISDEs 
define the same stochastic dynamics 
for sine$_2 $, Airy$_2 $, and Bessel$_2 $. 
\end{theorem}

By algebraic method, the finite particle approximation 
for sine$_2 $, Airy$_2 $, and Bessel$_2 $ is proved \cite{o-t.sm}. By analytic method, the same is proved for these point processes with 
$ \beta = 1,2,4$ and also the Ginibre point process \cite{k-o.sg,k-o.fpa}. The latter approach is robust and valid for many other examples.

\begin{acknowledgement}
H.O. is supported in part by a Grant-in-Aid for Scenic Research (KIBAN-A, No.\!\! 24244010; KIBAN-A, No.\!\! 16H02149; KIBAN-S, No.\!\! 16H06338) from the Japan Society for the Promotion of Science. 
\end{acknowledgement}

\input{is-ref.tex}
\end{document}

%% file: is-def.tex
\DN\lref[1]{Lemma~\ref{#1}}
\DN\tref[1]{Theorem~\ref{#1}}
\DN\pref[1]{Proposition~\ref{#1}}
\DN\sref[1]{Section~\ref{#1}}
\DN\dref[1]{Definition~\ref{#1}}
\DN\rref[1]{Remark~\ref{#1}} 
\DN\corref[1]{Corollary~\ref{#1}}
\DN\eref[1]{Example~\ref{#1}}
\DN\emp{\tc} 
 \DN\tcb[1]{#1}
 \DN\tcr[1]{#1}
 \DN\tcg[1]{#1}
 \DN\tcm[1]{}
 \DN\tcba[1]{}
\DN\bub{\tcb{$ \bullet $ }}
\DN\bur{\tcr{$ \bullet $ }}
\DN\bubb{\tcb{$ \bullet $ $ \bullet $ $ \bullet $ }}
\DN\burr{\tcr{$ \bullet $ $ \bullet $ $ \bullet $ }}

\DN\emr[1]{\emp{red}{#1}}
\DN\emg[1]{\emp{blue}{#1}}
\DN\emb[1]{\emp{blue}{#1}}

\DN\zzz{\end{slide}

%% file: is-5.bbl
\begin{thebibliography}{99}
\DN\ttl[1]{\textit{#1},} 

\bibitem{akr} Albeverio, S., Kondratiev, Yu. G., R\"{o}ckner M., 
\ttl{Analysis and geometry on configuration spaces:
the Gibbsian Case}, journal of functional analysis 157, 242-291 (1998)




\bibitem{Fr}Fritz, J., \ttl{Gradient Dynamics of Infinite Point Systems} Ann. Probab. {\bf 15} (1987) 478-514. 

\bibitem{FOT.2} Fukushima, M., {\it et al.}, \ttl{Dirichlet forms and symmetric Markov processes} 2nd ed., Walter de Gruyter (2011). 

\bibitem{o-h.bes} Honda,~R., Osada,~H, \ttl{Infinite-dimensional stochastic differential equations related to Bessel random point fields}, Stochastic Process. Appl.  {\bf 125}  (2015),  no. 10, 3801-3822.
 

\bibitem{GP} Ghosh,~S. and Peres,~Y, 
	Rigidity and Tolerance in point processes: Gaussian zeroes
	and Ginibre eigenvalues,  
	available at \url{http://arxiv.org/pdf/1211.2381v2.pdf}




\bibitem{j.02} Johansson, K., \ttl{Non-intersecting paths, random tilings and random matrices} Probab. Theory Relat. Fields {\bf 123}, 225-280 (2002).




\bibitem{KT07b} Katori, M., Tanemura, H.: \ttl{Noncolliding Brownian motion and determinantal processes}. J. Stat. Phys. {\bf 129}, 1233-1277 (2007). 


\bibitem{KT11}Katori, M., Tanemura, H., \ttl{Markov property of determinantal processes with extended sine, Airy, and Bessel kernels} Markov processes and related fields {\bf 17}, 541-580 (2011). 


\bibitem{KT11-b}Katori, M., Tanemura, H., \ttl{Noncolliding square Bessel processes} J.\ Stat. Phys. {\bf 142}, 592-615 (2011). 


\bibitem{k-o.fpa} Kawamoto, Y., Osada, H., \ttl{Finite-particle approximations for interacting Brownian particles with logarithmic potentials} (to appear in J. Math. Soc. Japan), arXiv:1607.06922

\bibitem{k-o.sg} Kawamoto, Y., Osada, H.: Dynamical bulk scaling limit of Gaussian unitary ensembles and stochastic-differential-equation gaps, (preprint),  arXiv:1610.05969


\bibitem{lang.1} Lang,~R., \ttl{Unendlich-dimensionale Wienerprocesse mit Wechselwirkung I}  Z. Wahrschverw. Gebiete  {\bf  38  } (1977) 55-72.   

\bibitem{lang.2} Lang,~R., \ttl{Unendlich-dimensionale Wienerprocesse mit Wechselwirkung II} Z. Wahrschverw. Gebiete  {\bf  39  } (1978) 277-299.   


\bibitem{m-r} Ma, Z.-M. and R\"ockner, M., \ttl{Introduction to the theory of (non-symmetric) Dirichlet forms} \ttl{Springer-Verlag} 1992. 



 

\bibitem{o.dfa} Osada, H., \ttl{Dirichlet form approach to infinite-dimensional Wiener processes with singular interactions} Commun. Math. Phys. {\bf 176}, 117-131 (1996). 



\bibitem{o.p} Osada,~H., \ttl{Positivity of the self-diffusion matrix of interacting Brownian particles with hard core} {\it Probab.\! Theory Relat.\!  Fields}, {\bf 112}, (1998), 53-90. 



\bibitem{o.tp} Osada,~H., \ttl{Tagged particle processes and their non-explosion criteria} J. Math. Soc. Japan, {\bf 62}, No.\ {\bf 3}, 867-894 (2010). 

\bibitem{o.isde} Osada,~H., \ttl{Infinite-dimensional stochastic differential equations related to random matrices} Probability Theory and Related Fields, {\bf 153}, 471-509 (2012). 

\bibitem{o.rm} Osada,~H., \ttl{Interacting Brownian motions in infinite dimensions with logarithmic interaction potentials} Ann. of  Probab. {\bf 41}, 1-49 (2013). 

\bibitem{o.rm2} Osada, H., \ttl{Interacting Brownian motions in infinite dimensions with logarithmic interaction potentials II : Airy random point field} Stochastic Processes and their applications {\bf 123}, 813-838 (2013). 


\bibitem{o.sub} Osada, H., \ttl{The Ginibre interacting Brownian motion is sub-diffusive} (preprint). 


\bibitem{o-o.tt} Osada, H., Osada, S., \ttl{Discrete approximations of determinantal point processes on continuous spaces: tree representations and  tail triviality}, arXive:1603.07478-v3.



\bibitem{o-s.abs} Osada,~H., Shirai,~T., \ttl{Absolute continuity and singularity of Palm measures of the Ginibre point process} Probab. Theory Relat. Fields, {\bf 165}, 725-770, DOI 10.1007/s00440-015-0644-6 



\bibitem{o-t.core} Osada,~H., Tanemura,~H., \ttl{Cores of Dirichlet forms related to Random Matrix Theory} Proc. Jpn. Acad., Ser. A, Vol. 90, 145-150 (2014). 

\bibitem{o-t.sm} Osada,~H., Tanemura,~H., \ttl{Strong Markov property of determinantal processes with extended kernels} Stochastic Processes and their Applications {\bf 126},  186-208, no. 1,  (2016). 

\bibitem{o-t.tail} Osada,~H., Tanemura,~H., \ttl{Infinite-dimensional stochastic differential equations and tail $ \sigma $-fields}, arXiv:1412.8674. 

\bibitem{o-t.airy} Osada,~H., Tanemura,~H., \ttl{Infinite-dimensional stochastic differential equations related to Airy random point fields}, arXiv:1408.0632.








\bibitem{shirai.06} Shirai, T., \textit{ Large deviations for the Fermion point process associated with the exponential kernel} J.\ Stat.\ Phys.\ {\bf 123}  (2006), 615-629. 





\bibitem{tane.2} Tanemura,~H., \ttl{A system of infinitely many mutually reflecting Brownian balls in $\R^d $}  Probab.\  Theory Relat.\  Fields {\bf  104} (1996) 399-426. 

\bibitem{tsai.14}  Tsai, Li-Cheng,  \ttl{Infinite dimensional stochastic differential equations for Dyson's model} Probab.\  Theory Relat.\  Fields  (published on line) DOI 10.1007/s00440-015-0672-2  (2015)


%
%


\bibitem{y.96} Yoshida, M.,W., \ttl{Construction of infinite-dimensional interacting diffusion processes through Dirichlet forms}, Probab. Theory Related Fields 106 (1996),265-297.

\end{thebibliography}
